\documentclass[12pt,a4paper]{article}
\usepackage{indentfirst}
\usepackage{amsmath,amssymb,amsfonts,amsthm,graphics}
\usepackage{makeidx}
\usepackage{amsmath,amssymb,amsfonts,amsthm,graphics,graphicx}
\usepackage{makeidx}
\usepackage{color}
\usepackage{ifpdf}
\usepackage{subfigure}
\usepackage{mathrsfs}

\newtheorem{thm}{Theorem}

\newtheorem{lem}{Lemma}
\newtheorem{pro}{Proposition}
\newtheorem{cor}{Corollary}

\theoremstyle{definition}

\newtheorem{example}{Example}

\def\-{\mbox{--}}

\newtheorem{remark}{Remark}

\def\pf{\noindent {\it Proof.} }
\begin{document}

\title{\Large\bf More on total monochromatic connection of graphs\footnote{Supported by NSFC No.11371205 and 11531011, and PCSIRT.} }
\author{\small Hui Jiang, Xueliang Li, Yingying Zhang\\
\small Center for Combinatorics and LPMC\\
\small Nankai University, Tianjin 300071, China\\
\small E-mail: jhuink@163.com; lxl@nankai.edu.cn;\\
\small zyydlwyx@163.com}
\date{}
\maketitle
\begin{abstract}

A graph is said to be {\it total-colored} if all the edges and the vertices of the graph are colored. A total-coloring of a graph is a {\it total monochromatically-connecting coloring} ({\it TMC-coloring}, for short) if any two vertices of the graph are connected by a path whose edges and internal vertices on the path have the same color. For a connected graph $G$, the {\it total monochromatic connection number}, denoted by $tmc(G)$, is defined as the maximum number of colors used in a TMC-coloring of $G$. Note that a TMC-coloring does not exist if $G$ is not connected, in which case we simply let $tmc(G)=0$. In this paper, we first characterize all graphs of order $n$ and size $m$ with $tmc(G)=3,4,5,6,m+n-2,m+n-3$ and $m+n-4$, respectively. Then we determine the threshold function for a random graph to have $tmc(G)\geq f(n)$, where $f(n)$ is a function satisfying $1\leq f(n)<\frac{1}{2}n(n-1)+n$. Finally, we show that for a given connected graph $G$, and a positive integer $L$ with $L\leq m+n$, it is NP-complete to decide whether $tmc(G)\geq L$.


\end{abstract}

\section{Introduction}

In this paper, all graphs are simple, finite and undirected. We refer to the book \cite{B} for undefined notation and terminology in graph theory. Throughout this paper, let $n$ and $m$ denote the order (number of vertices) and size (number of edges) of a graph, respectively. Moreover, a vertex of a connected graph is called a {\it leaf} if its degree is one; otherwise, it is an {\it internal vertex}. Let $l(T)$ and $q(T)$ denote the number of leaves and the number of internal vertices of a tree $T$, respectively, and let $l(G)=\max\{ l(T) | $ $T$ is a spanning tree of $G$ $\}$ and $q(G)=\min\{ q(T) | $ $T$ is a spanning tree of $G$ $\}$ for a connected graph $G$. Note that the sum of $l(G)$ and $q(G)$ is $n$ for any connected graph $G$ of order $n$. A path in an edge-colored graph is a {\it monochromatic path} if all the edges on the path have the same color. An edge-coloring of a connected graph is a {\it monochromatically-connecting coloring} ({\it MC-coloring}, for short) if any two vertices of the graph are connected by a monochromatic path of the graph. For a connected graph $G$, the {\it monochromatic connection number} of $G$, denoted by $mc(G)$, is defined as the maximum number of colors used in an MC-coloring of $G$. An {\it extremal MC-coloring} is an MC-coloring that uses $mc(G)$ colors. Note that $mc(G)=m$ if and only if $G$ is a complete graph. The concept of $mc(G)$ was first introduced by Caro and Yuster \cite{CY} and has been well-studied recently. We refer the reader to \cite{CLW,GLQZ} for more details.

In \cite{JLZ}, the authors introduced the concept of total monochromatic connection of graphs. A graph is said to be {\it total-colored} if all the edges and the vertices of the graph are colored. A path in a total-colored graph is a {\it total monochromatic path} if all the edges and internal vertices on the path have the same color. A total-coloring of a graph is a {\it total monochromatically-connecting coloring} ({\it TMC-coloring}, for short) if any two vertices of the graph are connected by a total monochromatic path of the graph. For a connected graph $G$, the {\it total monochromatic connection number}, denoted by $tmc(G)$, is defined as the maximum number of colors used in a TMC-coloring of $G$. Note that a TMC-coloring does not exist if $G$ is not connected, in which case we simply let $tmc(G)=0$. An {\it extremal TMC-coloring} is a TMC-coloring that uses $tmc(G)$ colors. It is easy to check that $tmc(G)=m+n$ if and only if $G$ is a complete graph. Actually, these concepts are not only inspired by the concept of monochromatic connection number but also by the concepts of monochromatic vertex connection number and total rainbow connection number of a connected graph. For details about them we refer to \cite{CLW1,JLZ1,LMS1,S}. From the definition of the total monochromatic connection number, the following results follow immediately.

\begin{pro}\label{pro1}\cite{JLZ} If $G$ is a connected graph and $H$ is a connected spanning subgraph of $G$, then $tmc(G)\geq e(G)-e(H)+tmc(H)$.
\end{pro}

\begin{thm}\label{thm1}\cite{JLZ} For a connected graph $G$, $tmc(G)\geq m-n+2+l(G)$.
\end{thm}

In particular, $tmc(G)=m-n+2+l(G)$ if $G$ is a tree. The authors \cite{JLZ} also showed that there are dense graphs that still meet this lower bound.

\begin{thm}\label{thm2}\cite{JLZ} Let $G$ be a connected graph of order $n>3$. If $G$ satisfies any of the following properties, then $tmc(G)=m-n+2+l(G)$.

$(a)$ The complement $\overline{G}$ of $G$ is $4$-connected.

$(b)$ $G$ is $K_3$-free.

$(c)$ $\Delta(G)<n-\frac{2m-3(n-1)}{n-3}$.

$(d)$ $diam(G)\geq 3$.

$(e)$ $G$ has a cut vertex.
\end{thm}

Moreover, the authors \cite{JLZ} gave an example to show that the lower bound $m-n+2+l(G)$ is not always attained.

\begin{example}\label{example1}\cite{JLZ} Let $G= K_{n_1,\ldots,n_r}$ be a complete multipartite graph with $n_1 \geq \ldots \geq n_t\geq 2$ and $n_{t+1}=\ldots=n_r=1$. Then $tmc(G)=m+r-t$.
\end{example}

Let $G$ be a connected graph and $f$ be an extremal TMC-coloring of $G$ that uses a given color $c$. Note that the subgraph $H$ formed by the edges and vertices with color $c$ is a tree where the color of each internal vertex is $c$ \cite{JLZ}. Now we define the {\it color tree} as the tree formed by the edges and vertices with color $c$, denoted by $T_c$. If $T_c$ has at least two edges, the color $c$ is called {\it nontrivial}; otherwise, $c$ is {\it trivial}. We call an extremal TMC-coloring {\it simple} if for any two nontrivial colors $c$ and $d$, the corresponding trees $T_c$ and $T_d$ intersect in at most one vertex. If $f$ is simple, then the leaves of $T_c$ must have distinct colors different from color $c$. Moreover, a nontrivial color tree of $f$ with $m'$ edges and $q'$ internal vertices is said to {\it waste} $m'-1+q'$ colors. For the rest of this paper we will use these facts without further mentioning them. In addition, we list a helpful lemma below.

\begin{lem}\label{lem1}\cite{JLZ} Every connected graph $G$ has a simple extremal TMC-coloring.
\end{lem}

This paper is organized as follows. In Section $2$, we characterize all graphs $G$ with $tmc(G)=3,4,5,6,m+n-2,m+n-3,m+n-4$, respectively. In Section $3$, we show that for any function $f(n)$ satisfying $1\leq f(n)<\frac{1}{2}n(n-1)+n$, if $ln\log n\leq f(n)< \frac{1}{2}n(n-1)+n$, where $l\in\mathbb{R}^{+}$, then $p=\frac{f(n)+n\log\log n}{n^2}$ is a sharp threshold function for the property $tmc(G(n,p))\geq f(n)$; if $f(n)=o(n\log n)$, then $p=\frac{\log n}{n}$ is a sharp threshold function for the property $tmc(G(n,p))\geq f(n)$. In Section $4$, we prove that for a given connected graph $G$, and a positive integer $L$ with $L\leq m+n$, it is NP-complete to decide whether $tmc(G)\geq L$.

\section{Characterization of graphs with small or large $tmc$}

In this section, we characterize all graphs $G$ with $tmc(G)=3,4,5,6,m+n-2,m+n-3,m+n-4$, respectively. We call a connected graph $G$ {\it unicyclic}, {\it bicyclic}, or {\it tricyclic} if $m=n,n+1$ or $n+2$, respectively. Let $\mathcal{T}_i$ denote the set of the trees with $l(G)=i$, where $2\leq i\leq n-1$. Note that if $G$ is a connected graph with $l(G)=2$, then $G$ is either a path or a cycle.

\begin{thm}\label{thm3} Let $G$ be a connected graph. Then $tmc(G)=3$ if and only if $G$ is a path.
\end{thm}

\pf If $G$ is a path, then $tmc(G)=m-n+2+l(G)=3$. Hence it remains to verify the
converse. Let $G$ be a connected graph with $tmc(G)=3$. By Theorem \ref{thm1}, we get that $m\leq n+1-l(G)$ and then $m\leq n-1$ as $l(G)\geq2$. Since $G$ is a connected graph, it follows that $m=n-1$ and $l(G)=2$. Thus $G$ is a path.
\qed

\begin{thm}\label{thm4} Let $G$ be a connected graph. Then $tmc(G)=4$ if and only if $G\in\mathcal{T}_3$ or $G$ is a cycle except for $K_3$.
\end{thm}

\pf If $G\in\mathcal{T}_3$ or $G$ is a cycle except for $K_3$, then $tmc(G)=4$ by Theorem \ref{thm2}$(b)$. Conversely, let $G$ be a connected graph with $tmc(G)=4$. First, we have $m\leq n+2-l(G)$ by Theorem \ref{thm1}. Since $l(G)\geq2$ and $m\geq n-1$, it follows that $l(G)=2$ or $3$. If $l(G)=3$, then $m=n-1$ and so $G\in\mathcal{T}_3$. Otherwise, from Theorem \ref{thm3} we have that $G$ is a cycle and $G\neq K_3$ since $tmc(K_3)=6$. \qed

\begin{figure}[h,t,b,p]
\begin{center}
\scalebox{0.8}[0.8]{\includegraphics{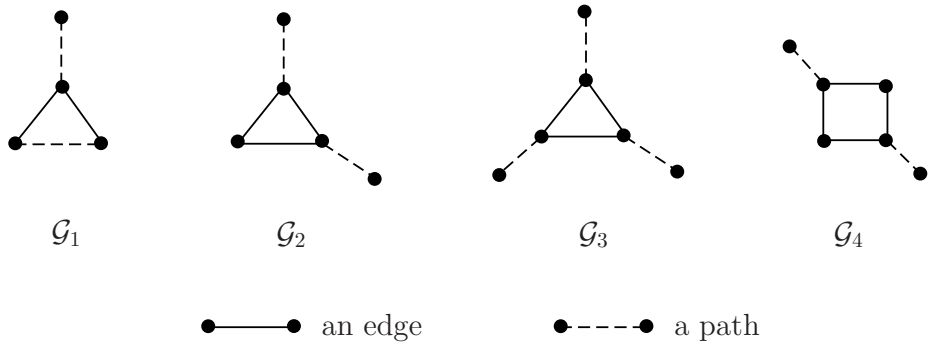}}
\end{center}
\caption{Unicyclic graphs with $l(G)=3$.}\label{Fig.1.}
\end{figure}

\begin{thm}\label{thm5} Let $G$ be a connected graph. Then $tmc(G)=5$ if and only if $G\in\mathcal{T}_4$ or $G\in\mathcal{G}_i$, where $1\leq i\leq4$; see Figure \ref{Fig.1.}.
\end{thm}

\pf If $G\in\mathcal{T}_4$ or $G\in\mathcal{G}_i$, where $1\leq i\leq4$, then $G$ has a cut vertex and so $tmc(G)=5$ by Theorem \ref{thm2}$(e)$. Hence it remains to verify the
converse. Let $G$ be a connected graph with $tmc(G)=5$. First, we have $m\leq n+3-l(G)$ by Theorem \ref{thm1}. Since $l(G)\geq2$ and $m\geq n-1$, it follows that $l(G)=2,3$ or $4$. If $l(G)=4$, then $m=n-1$ and so $G\in\mathcal{T}_4$. If $l(G)=3$, then we have $m=n$ from Theorem \ref{thm4} and so $G$ is a unicyclic graph with $l(G)=3$; see Figure \ref{Fig.1.}. If $l(G)=2$, then we have $G=K_3$ from Theorems \ref{thm3} and \ref{thm4}. However, $tmc(K_3)=6$, a contradiction.\qed

\begin{figure}[h,t,b,p]
\begin{center}
\scalebox{0.8}[0.8]{\includegraphics{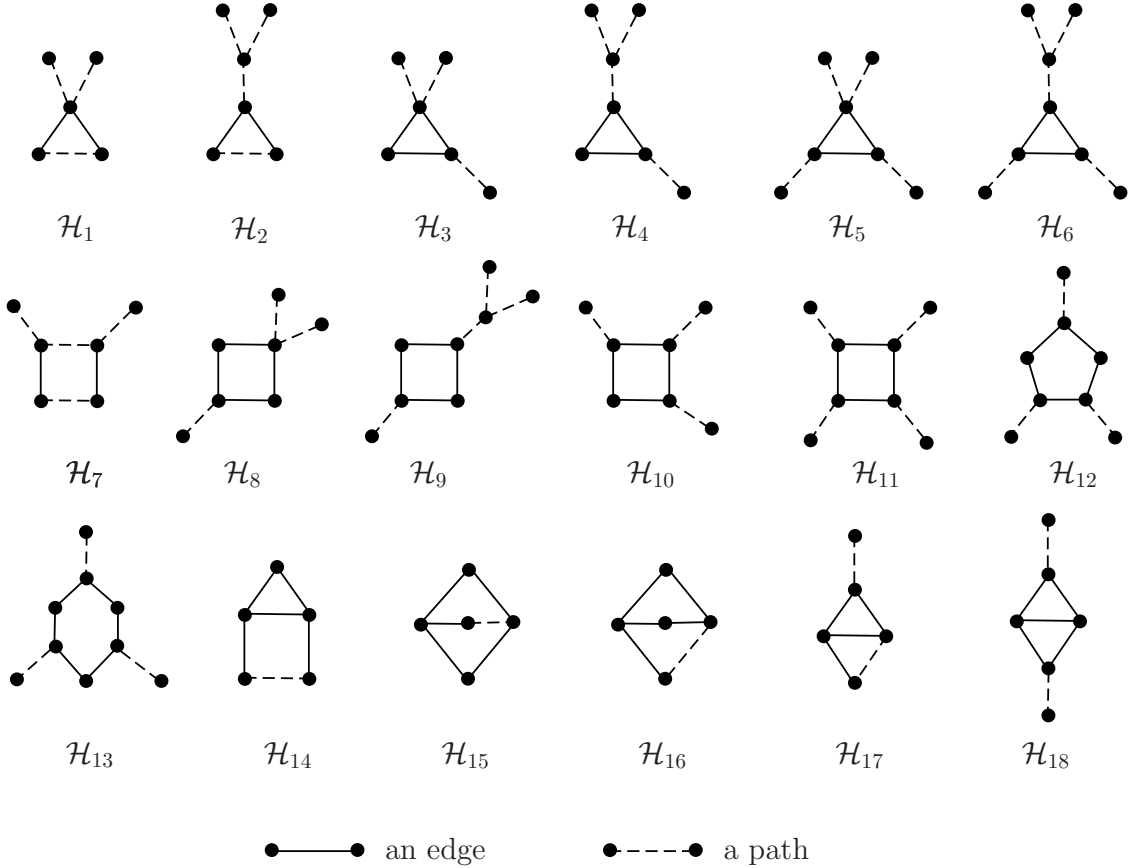}}
\end{center}
\caption{The graphs in Theorem 6.}\label{Fig.2.}
\end{figure}

\begin{thm}\label{thm6} Let $G$ be a connected graph. Then $tmc(G)=6$ if and only if $G=K_3$, $G\in\mathcal{T}_5$ or $G\in\mathcal{H}_i$, where $1\leq i\leq18$; see Figure \ref{Fig.2.}.
\end{thm}

\pf It is easy to verify the sufficiency by Theorem \ref{thm2}. Next we just need to prove the necessity. Let $G$ be a connected graph with $tmc(G)=6$. First, we have $m\leq n+4-l(G)$ by Theorem \ref{thm1}. Since $l(G)\geq2$ and $m\geq n-1$, it follows that $l(G)=2,3,4$ or $5$. If $l(G)=5$, then $m=n-1$ and so $G\in\mathcal{T}_5$. If $l(G)=4$, we have that $m=n$ from Theorem \ref{thm5} and so $G$ is a unicyclic graph with $l(G)=4$; see $\mathcal{H}_i\ (1\leq i\leq13)$ in Figure \ref{Fig.2.}. Similarly, from Theorems \ref{thm4} and \ref{thm5}, we have that $m=n+1$ if $l(G)=3$ and then $G$ is a bicyclic graph with $l(G)=3$ except for $K_{2,1,1}$ since $tmc(K_{2,1,1})=7$; see $\mathcal{H}_i\ (14\leq i\leq18)$ in Figure \ref{Fig.2.}. If $l(G)=2$, we have that $G=K_3$ from Theorems \ref{thm3} and \ref{thm4}. \qed

Recall that $tmc(G)=m+n$ if and only if $G=K_n$. In fact, there does not exist a graph such that $tmc(G)=m+n-1$. We are given a connected graph $G$ with $diam(G)\geq2$ and a simple extremal TMC-coloring of $G$. Then there must be two nonadjacent vertices in a nontrivial color tree. Since every nontrivial color tree wastes at least 2 colors, we get that $tmc(G)\leq m+n-2$. Hence in the following, we characterize all graphs $G$ having $tmc(G)=m+n-2,m+n-3,m+n-4$. Let $tK_2$ be $t$ nonadjacent edges of $K_n$, where $t\leq\lfloor\frac{n}{2}\rfloor$. Given a graph $H$, let $K_n-H$ denote the graph obtained from $K_n$ by deleting the edges of $H$.

\begin{thm}\label{thm7} Let $G$ be a connected graph. Then $tmc(G)=m+n-2$ if and only if $G=K_n-K_2$.
\end{thm}

\pf Clearly, $K_n-K_2=K_{2,1,...,1}$. Then $tmc(K_n-K_2)=m+n-2$. Conversely, let $G$ be a connected graph with $tmc(G)=m+n-2$. We are given a simple extremal TMC-coloring $f$ of $G$. Suppose that $f$ consists of $k$ nontrivial color trees, denoted by $T_1,\ldots,T_k$. Since each nontrivial color tree wastes at least two colors, it follows that $k=1$ and $T_1=P_3$. Thus, $G=K_n-K_2$.\qed

\begin{thm}\label{thm8} Let $G$ be a connected graph. Then $tmc(G)=m+n-3$ if and only if $G$ is either $K_n-K_3$ or $K_n-P_3$.
\end{thm}

\pf Note that $K_n-K_3=K_{3,1,...,1}$ and then $tmc(K_n-K_3)=m+n-3$. Note that $K_n-K_3$ is a spanning subgraph of $K_n-P_3$. Then $tmc(K_n-P_3)\geq m+n-3$ by Proposition \ref{pro1}. Now we just need to prove that $tmc(K_n-P_3)\leq m+n-3$. Let $f$ be a simple extremal TMC-coloring of $K_n-P_3$. Suppose that $f$ consists of $k$ nontrivial color trees. Since there are two pairs of nonadjacent vertices in two nontrivial color trees or in a common nontrivial color tree, it wastes at least three colors and then $tmc(K_n-P_3)\leq m+n-3$. Hence $tmc(K_n-P_3)=m+n-3$.

Now it remains to verify the converse. Let $G$ be a connected graph with $tmc(G)=m+n-3$. We are given a simple extremal TMC-coloring $f$ of $G$. Suppose that $f$ consists of $k$ nontrivial color trees, denoted by $T_1,\ldots,T_k$. Since each nontrivial color tree wastes at least two colors, we get that $k=1$ and $T_1=K_{1,3}$. Thus, $K_n-K_3$ is a spanning subgraph of $G$. From Theorem \ref{thm7}, it can be checked that $G$ is either $K_n-K_3$ or $K_n-P_3$.
\qed

\begin{thm}\label{thm9} Let $G$ be a connected graph. Then $tmc(G)=m+n-4$ if and only if $G\in\{K_n-P_4,K_n-2K_2,K_n-K_4,K_n-(K_4-K_2),K_n-(K_4-P_3),K_n-C_4,K_n-K_{1,3}\}$.
\end{thm}

\pf Clearly, $K_n-2K_2=K_{2,2,1,...,1}$ and $K_n-K_4=K_{4,1,...,1}$. Thus we have that $tmc(K_n-2K_2)=tmc(K_n-K_4)=m+n-4$. If $G=K_n-P_4$, there are three pairs of nonadjacent vertices and let $f$ be a simple extremal TMC-coloring of $G$. Suppose that $f$ consists of $k$ nontrivial color trees. Then it wastes at least 4 colors and so $tmc(G)\leq m+n-4$. Since $K_n-K_4$ is a spanning subgraph of $G$, $tmc(G)\geq m+n-4$ by Proposition \ref{pro1}. Thus we get that $tmc(K_n-P_4)=m+n-4$. Similarly, it can be verified that $tmc(K_n-(K_4-K_2))=tmc(K_n-(K_4-P_3))=tmc(K_n-C_4)=tmc(K_n-K_{1,3})=m+n-4$.

Conversely, let $G$ be a connected graph with $tmc(G)=m+n-4$. We are given a simple extremal TMC-coloring $f$ of $G$. Suppose that $f$ consists of $k$ nontrivial color trees, denoted by $T_1,\ldots,T_k$. Since each nontrivial color tree wastes at least two colors, we get the following two cases.

\textbf{Case 1.} $k=1$.

Then $T_1=P_4$ or $K_{1,4}$. If $T_1=P_4$, then $K_n-P_4$ is a spanning subgraph of $G$. From Theorems \ref{thm7} and \ref{thm8}, we obtain that $G$ is either $K_n-P_4$ or $K_n-2K_2$. If $T_1=K_{1,4}$, then $K_n-K_4$ is a spanning subgraph of $G$. From Theorems \ref{thm7} and \ref{thm8}, we get that $G\in\{K_n-P_4,K_n-2K_2,K_n-K_4,K_n-(K_4-K_2),K_n-(K_4-P_3),K_n-C_4,K_n-K_{1,3}\}$.

\textbf{Case 2.} $k=2$.

Then $T_1=T_2=P_3$. From Theorem \ref{thm8}, $T_1$ and $T_2$ have not a common leaf. Thus $K_n-2K_2$ is a spanning subgraph of $G$. Since $tmc(K_n-K_2)=m+n-2$ and $tmc(K_n)=m+n$, we have that $G=K_n-2K_2=K_{2,2,1,...,1}$.
\qed

\section{Random graphs}

Let $G=G(n,p)$ denote the random graph with $n$ vertices and edge probability $p$ \cite{AS}. For a graph property $P$ and for a function $p=p(n)$, we say that $G(n,p)$ satisfies $P$ {\it almost surely} if the probability that $G(n,p(n))$ satisfies $P$ tends to $1$ as $n$ tends to infinity. We say that a function $f(n)$ is a {\it sharp threshold function} for the property $P$ if there are two positive constants $C$ and $c$ so that $G(n,p)$ satisfies $P$ almost surely for all $p\geq Cf(n)$ and $G(n,cf(n))$ almost surely does not satisfy $P$.

Let $G$ and $H$ be two graphs on $n$ vertices. A property $P$ is said to be {\it monotone} if whenever $G\subseteq H$ and $G$ satisfies $P$, then $H$ also satisfies $P$. It is well-known that all monotone graph properties have sharp threshold functions; see
\cite{BT} and \cite{FK}. For any graph $G$ with $n$ vertices and any function $f(n)$, having $tmc(G)\geq f(n)$ is a monotone graph property (adding edges does not destroy this property), so it has a sharp threshold function. In the following, we establish a sharp threshold function for the graph property $tmc(G)\geq f(n)$.

\begin{thm}\label{thm10} Let $f(n)$ be a function satisfying $1\leq f(n)<\frac{1}{2}n(n-1)+n$. Then

\begin{eqnarray}p=
\begin{cases}
\frac{f(n)+n\log\log n}{n^2} &if\ ln\log n\leq f(n)< \frac{1}{2}n(n-1)+n,\ where\ l\in \mathbb{R}^{+}, \cr
\frac{\log n}{n} &if\ f(n)=o(n\log n).\end{cases}
\end{eqnarray}
is a sharp threshold function for the property $tmc(G(n,p))\geq f(n)$.
\end{thm}

\begin{remark}\label{rem1} Note that if $f(n)=\frac{1}{2}n(n-1)+n$, then $G(n,p)$ is a complete graph $K_n$ and $p=1$. Hence we only concentrate on the case $f(n)<\frac{1}{2}n(n-1)+n$.
\end{remark}

Before proving Theorem \ref{thm10}, we need some lemmas.

\begin{lem}\label{lem2}\cite{ER} Let $p=\frac{\log n+a}{n}$. Then

\begin{eqnarray} Pr[G(n,p)\ is\ connected]\rightarrow
\begin{cases}
e^{e^{-a}} &if\ |a|=O(1), \cr
0 &a\rightarrow -\infty, \cr
1 &a\rightarrow +\infty.
\end{cases}
\end{eqnarray}
\end{lem}

\begin{lem}\label{lem3} \cite{AS} {\bf (Chernoff Bound)} If $X$ is binomial variable with expectation $\mu$, and $0<\delta<1$, then
$$Pr[X<(1-\delta)\mu)]\leq exp(-\frac{\delta^2\mu}{2})$$
and
$$Pr[X>(1+\delta)\mu)]\leq exp(-\frac{\delta^2\mu}{2+\delta}).$$
\end{lem}

\begin{lem}\label{lem4}  Let $G$ be a noncomplete connected graph of order $n$ with minimum degree $\delta$. Then $tmc(G)\leq m-n+\delta+1+l(G)$.
\end{lem}

\pf For a noncomplete graph $G$, we have that $tmc(G)\leq mc(G)+l(G)$ whose proof is contained in the proof of Theorem 6 in \cite{JLZ}. Moreover, $mc(G)\leq m-n+\delta+1$ by Proposition 12 in \cite{CY}. Thus $tmc(G)\leq m-n+\delta+1+l(G)$.\qed

\noindent {\bf Proof of Theorem 10:} We divide our proof into two cases according to the range of $f(n)$.

\textbf{Case 1.}  $ln\log n\leq f(n)< \frac{1}{2}n(n-1)+n$, where $l\in \mathbb{R}^{+}$.

We first prove that there exists a constant $C$ such that the random graph $G(n,Cp)$ with $p=\frac{f(n)+n\log\log n}{n^2}$ almost surely has $tmc(G(n,Cp))\geq f(n)$. Let
\begin{eqnarray} C=
\begin{cases}
5 &if\ n\log n\leq f(n)< \frac{1}{2}n(n-1)+n, \cr
\frac{5}{l} &if\ f(n)=ln\log n,\ where\ 0<l<1.
\end{cases}
\end{eqnarray} It is easy to check that $G(n,Cp)$ is almost surely connected by Lemma \ref{lem2}. Let $\mu_1$ denote the expectation of the number of edges in $G(n,Cp)$. Then
$$\mu_1=\frac{n(n-1)}{2}\cdot Cp=\frac{C}{2}(\frac{n-1}{n}f(n)+(n-1)\log\log n).$$
Moreover from Lemma \ref{lem3}, it follows that $Pr[|E(G(n,Cp))|<\frac{\mu_1}{2}]\leq\exp(-\frac{1}{8}\mu_1)=o(1)$. Suppose that $|E(G(n,Cp))|\geq \frac{\mu_1}{2}$. By Theorem \ref{thm1} we have that for $n$ sufficiently large,

$$tmc(G(n,Cp))\geq |E(G(n,Cp))|-n+2+l(G)\geq \frac{\mu_1}{2}-n+2+l(G)$$
$$\ \ \ \ \ \ \ \ \ \ \ \ \ \ \ \ \ =\frac{C}{4}(\frac{n-1}{n}f(n)+(n-1)\log\log n)-n+2+l(G)$$
$$\ \ \ \ \ \ \ \ \ \ \ \ \geq \frac{5}{4}(\frac{n-1}{n}f(n)+(n-1)\log\log n)-n+2+2$$
$$\ \ \ \ \ \ \ \geq f(n).\ \ \ \ \ \ \ \ \ \ \ \ \ \ \ \ \ \ \ \ \ \ \ \ \ \ \ \ \ \ \ \ \ \ \ \ \ \ \ \ \ \ \ \ \ \ $$
Thus, we conclude that $tmc(G(n,Cp))\geq f(n)$ holds with the probability at least $1-\exp(-\frac{1}{8}\mu_1)=1-o(1)$.

Next we show that there exists a constant $c$ such that the random graph $G(n,cp)$ with $p=\frac{f(n)+n\log\log n}{n^2}$ almost surely has $tmc(G(n,cp))<f(n)$. Let $c=1$ and $\mu_2$ denote the expectation of the number of edges in $G(n,cp)$. Then we have
$$\mu_2=\frac{n(n-1)}{2}\cdot cp=\frac{1}{2}(\frac{n-1}{n}f(n)+(n-1)\log\log n).$$
Furthermore by Lemma \ref{lem3}, it follows that $Pr[|E(G(n,cp))|>\frac{3}{2}\mu_2]\leq\exp(-\frac{1}{10}\mu_2)=o(1)$.
If $G(n,cp)$ is not connected, then $tmc(G(n,cp))=0<f(n)$. Otherwise, let $\delta$ be the minimum degree of $G(n,cp)$. Suppose that $|E(G(n,p))|\leq \frac{3}{2}\mu_2$. From Lemma \ref{lem4}, we have that for $n$ sufficiently large,
$$tmc(G(n,cp))\leq |E(G(n,cp))|-n+\delta+1+l(G)\leq \frac{3}{2}\mu_2-n+\delta+1+l(G)$$
$$\ \ \ \ \ \ \ \ \ \  =\frac{3}{4}(\frac{n-1}{n}f(n)+(n-1)\log\log n)-n+\delta+1+l(G)$$
$$\ \ \ \ \ \ \ \ \ \ \ \ < \frac{3}{4}(\frac{n-1}{n}f(n)+(n-1)\log\log n)-n+n+1+n-1$$
$$\ \ \ \ \ \ \ < f(n).\ \ \ \ \ \ \ \ \ \ \ \ \ \ \ \ \ \ \ \ \ \ \ \ \ \ \ \ \ \ \ \ \ \ \ \ \ \ \ \ \ \ \ \ \ \ \ \ \ \ \ \ \ \ \ \ \ $$
Hence, we conclude that $tmc(G(n,cp))< f(n)$ holds with the probability at least $1-\exp(-\frac{1}{10}\mu_2)=1-o(1)$.

\textbf{Case 2.} $f(n)=o(n\log n)$.

Let $C=2$ and $c=\frac{1}{2}$. By Lemma \ref{lem2}, we have that $G(n,Cp)$ is almost surely connected and $G(n,cp)$ is almost surely not connected. It can be checked that $tmc(G(n,Cp))\geq f(n)$ almost surely holds in a similar way as Case 1. On the other hand, since $G(n,cp)$ is almost surely not connected, $tmc(G(n,cp))=0<f(n)$ almost surely holds.
\qed

\section{Hardness result for computing $tmc$}

Given a graph $G$, a set $D\subseteq V(G)$ is called a {\it dominating set} of $G$ if every vertex of $G$ not in $D$ has a neighbor in $D$. If the subgraph induced by $D$ is connected, then $D$ is called a {\it connected dominating set}. The {\it connected dominating number}, denoted by $\gamma_c(G)$, is the minimum cardinalities of the connected dominating sets of $G$. Note that the sum of $\gamma_c(G)$ and $l(G)$ is $n$ because a vertex subset is a connected dominating set if and only if its complement is contained in the set of leaves of a spanning tree. In this section, we mainly prove the following result.

\begin{thm}\label{thm11} The following problem is NP-complete: Given a connected graph $G$ and a positive integer $L\leq m+n$, decide whether $tmc(G)\geq L$.
\end{thm}

In order to prove Theorem \ref{thm11}, we need the lemma as follows.

\begin{lem}\label{lem5}\cite{LW} The first problem defined below is polynomially reducible to the second one:
\\ Problem 1. Given a graph $G$ and a positive integer $k\leq n$, decide whether there is a dominating set of size $k$ or less.
\\ Problem 2. Given a connected graph $G$ with a cut vertex and a positive integer $k$ with $k\leq n$, decide whether there is a connected dominating set of size $k$ or less.
\end{lem}

\noindent {\bf Proof of Theorem 11:} Given a connected graph $G$ with a cut vertex, and a positive integer $k\leq n$. Note that $\gamma_c(G)\leq k$ if and only if $tmc(G)=m-n+2+l(G)=m-\gamma_c(G)+2\geq m-k+2$ by Theorem \ref{thm2}$(e)$. Then Problem 2 can be polynomially reducible to Problem 3: given a connected graph $G$ with a cut vertex and a positive integer $L$ with $L\leq m+n$, decide whether $tmc(G)\geq L$. Thus, Problem 1 can be reducible to Problem 3 by Lemma \ref{lem5}. Moreover, Problem 1 is known as a NP-complete problem in \cite{GJ}. Hence the problem in Theorem \ref{thm11} is NP-hard.

Next we prove that given a connected graph $G$ and a nonnegative integer $K<m+n$, to decide whether $tmc(G)\geq m+n-K$ is NP. Recall that a problem belongs to NP-class if given any instance of the problem whose answer is ``yes", there is a certificate validating this fact which can be checked in polynomial time. For any fixed integer $K$, to prove the problem of deciding whether $tmc(G)\geq m+n-K$ is NP,
we choose a TMC-coloring of $G$ with $m+n-K$ colors as a certificate. For checking a TMC-coloring with $m+n-K$ colors, we only need to check that $m+n-K$ colors are used and for any two vertices $u$ and $v$ of $G$, there exists a total monochromatic path between them. Notice that for any two vertices $u$ and $v$ of $G$, there are at most $n^{l-1}$ paths of length $l$, since if we let $P=uv_1v_2\cdots v_{l-1}v$, then there are less than $n$ choices for each $v_i\ (i\in\{1,2,\ldots,l-1\})$. Clearly, the path $P$ wastes at least $2l-2$ colors. Then $tmc(G)\leq m+n-(2l-2)$ and so $m+n-K\leq m+n-(2l-2)$ which implies that $l\leq\frac{K+2}{2}$. Therefore, $G$ contains at most
$\sum_{l=1}^{\frac{K+2}{2}} n^{l-1}=\frac{n^{\frac{K+2}{2}}-1}{n-1}=O(n^{\frac{K}{2}})$ $u$-$v$ paths of length at most $\frac{K+2}{2}$. Then, check these paths in turn until one finds a path whose edges and internal vertices have the same color. It follows that the time used for checking is at most $O(n^{\frac{K}{2}}\cdot n^2\cdot n^2)=O(n^{\frac{K}{2}+4})$. Since $K$ is a fixed integer, we conclude that the certificate can be checked in polynomial time. Then the problem of deciding whether
$tmc(G)\geq m+n-K$ belongs to NP-class and so is the problem in Theorem \ref{thm11}.

Therefore, the proof is complete.
\qed

\begin{cor}\label{cor1} Let $G$ be a connected graph. Then computing $tmc(G)$ is NP-hard.
\end{cor}

\end{document}